# Exact alternative optima for nonlinear optimization problems defined with maximum component objective function constrained by the Sugeno-Weber fuzzy relational inequalities

**Amin Ghodousian[1], Sara Zal[2], Minoo Ahmadi[3]**
A School of Engineering Science, College of Engineering University of Tehran,
Tehran, Iran

**Abstract**
In this paper, we study a latticized optimization problem with fuzzy relational inequality constraints where the feasible region is formed as the intersection of two inequality fuzzy systems and Sugeno-Weber family of t-norms is considered as fuzzy composition. Sugeno-Weber family of t-norms and t-conorms is one of the most applied one in various fuzzy modelling problems. This family of t-norms and t-conorms was suggested by Weber for modeling intersection and union of fuzzy sets. Also, the t-conorms were suggested as addition rules by Sugeno for so-called –fuzzy measures. The resolution of the feasible region of the problem is firstly investigated when it is defined with max-Sugeno-Weber composition and a necessary and sufficient condition is presented for determining the feasibility. Then, based on some theoretical properties of the problem, an algorithm is presented for solving this nonlinear problem. It is proved that the algorithm can find the exact optimal solution and an example is presented to illustrate the proposed algorithm.

# 1. Introduction

Resolution of fuzzy relational equations (FRE) with max-min composition was first studied by Sanchez [16]. Besides, Sanchez developed the application of FRE in medical diagnosis in biotechnology. Nowadays, it is well known that many issues associated with a body knowledge can be treated as FRE problems [14]. Since then the composition operator in FREs was replaced by max-product and furthermore extended to the general max-t-norm composition operator [2,5,6,8,9]. The resolution method was kept in improving. In fact, if the max-t-norm fuzzy relational equations is consistent, then its solution set is often a non-convex set that is completely determined by a unique maximum solution and a finite number of minimal solutions [3,8,9,15].

Most of the existing literatures in this research field focused on the resolution of FRE and its relevant optimization problems. However, only a few research works investigated the fuzzy relational inequalities (FRI) and its relevant optimization problems [5-7,10,11,18]. For example, Guo et al. [10] studied the linear programming problem with max-min FRI constraint. Based on the concept of pseudo-minimal index, Yang [17] developed a pseudo-minimal-index algorithm to minimize a linear objective function with addition-min FRI constraint. To improve the results presented in [17], Yang et al. [19] proposed the min-max programming subject to addition-min fuzzy relational inequalities. They also studied the multi-level linear programming problem with addition-min FRI constraint [20].

The optimization problems with general nonlinear objective functions and FRE or FRI constraints were studied in [4-9]. In general, the genetic algorithm was applied to deal with this kind of problems. However, some fuzzy relation nonlinear optimization problems could be solved by some specific method. For example, fuzzy relation geometric programming problem was introduced by Yang and Cao [1]. Yang et al. [21] studied the single-variable term semi-latticized geometric programming subject to max-product fuzzy relation equations. The proposed problem was devised from the peer-to-peer network system and the target was to minimize the biggest dissatisfaction degrees of the terminals in such system. Yang et al. [18] introduced another version of the latticized programming problem subject to max-prod fuzzy relation inequalities with application in the optimization management model of wireless communication emission base stations.

The latticized problem was defined by minimizing objective function $z(x) = \max\{x_1, x_2, \ldots, x_n\}$ subject to feasible region $X(A, b) = \{x \in [0,1]^n : A \circ x \geq b\}$ where "$\circ$" denotes fuzzy max-product composition. In this paper, we study the following non-linear optimization problem in which the objective function is defined as the maximum components function and the constraints are formed as as the intersection of two fuzzy systems of relational inequalities defined by Sugeno-Weber family of t-norms:

$$\begin{aligned} & min \ z(x) = \max\{x_1, x_2, \ldots, x_n\} \\ & A\varphi x \leq b^1 \\ & D\varphi x \geq b^2 \\ & x \in [0,1]^n \end{aligned} \qquad (1)$$

where $I_1 = \{1, 2, \ldots, m_1\}$, $I_2 = \{m_1 + 1, m_1 + 2, \ldots, m_1 + m_2\}$ and $J = \{1, 2, \ldots, n\}$. $A = (a_{ij})_{m_1 \times n}$ and $D = (d_{ij})_{m_2 \times n}$ are fuzzy matrices such that $0 \leq a_{ij} \leq 1$ ($\forall i \in I_1$ and $\forall j \in J$) and $0 \leq d_{ij} \leq 1$ ($\forall i \in I_2$ and $\forall j \in J$). $b^1 = (b_i^1)_{m_1 \times 1}$ is an $m_1$–dimensional fuzzy vector in $[0,1]^{m_1}$ (i.e., $0 \leq b_i^1 \leq 1$, $\forall i \in I_1$), $b^2 = (b_i^2)_{m_2 \times 1}$ is an $m_2$–dimensional fuzzy vector in $[0,1]^{m_2}$ (i.e., $0 \leq b_i^2 \leq 1$, $\forall i \in I_2$) and "$\varphi$" is the max-Sugeno-Weber composition, that is,

$$\varphi(x, y) = T_{SW}^{\lambda}(x, y) = \max\left\{\frac{x + y - 1 + \lambda xy}{1 + \lambda}, 0\right\}$$

in which $\lambda > -1$.

By these notations, problem (1) can be also expressed as follows

$$\begin{aligned} & min \ z(x) = \max\{x_1, x_2, \ldots, x_n\} \\ & \max_{j \in J}\{T_{SW}^{\lambda}(a_{ij}, x_j)\} \leq b_i^1 \quad, i \in I_1 \\ & \max_{j \in J}\{T_{SW}^{\lambda}(d_{ij}, x_j)\} \geq b_i^2 \quad, i \in I_2 \\ & x \in [0,1]^n \end{aligned} \qquad (2)$$

The members of the family $\{T_{SW}^{\lambda}\}$ are increasing functions of the parameter $\lambda$. It can be easily shown that Sugeno-Weber t-norm $T_{SW}^{\lambda}(x, y)$ converges to the product fuzzy intersection $xy$ as $\lambda$ goes to infinity and converges to Drastic product t-norm as $\lambda$ approaches $-1$.

The rest of the paper is organized as follows. In Section 2, a necessary and sufficient condition is derived to determine the feasibility of max-Sugeno-Weber fuzzy relational inequalities (FRI). In Section 3, the feasible solution set of Problem (1) is characterized. It is shown that the feasible region can be expressed as the union of a finite closed convex cells. Section 4 describes the optimal solution of Problem (1). An algorithm is proposed to find the optimal solution and finally, Section 5 provides a numerical example to illustrate the algorithm.

**2. Basic properties of max-Sugeno-Weber FRI**

Let $S_{T_{SW}^\lambda}(A, b^1)$ and $S_{T_{SW}^\lambda}(D, b^2)$ denote the feasible solutions sets of inequalities $A\varphi x \leq b^1$ and $D\varphi x \geq b^2$, respectively. Also, let $S_{T_{SW}^\lambda}(A, D, b^1, b^2)$ denote the feasible solutions set of problem (1). Based on the foregoing notations, it is clear that $S_{T_{SW}^\lambda}(A, D, b^1, b^2) = S_{T_{SW}^\lambda}(A, b^1) \cap S_{T_{SW}^\lambda}(D, b^2)$.

**Definition 1.** Define $J_i^1 = \{j \in J : a_{ij} \geq b_i^1\}$, $\forall i \in I_1$, and $J_i^2 = \{j \in J : d_{ij} \geq b_i^2\}$, $\forall i \in I_2$.

**Definition 2.** For each $i \in I_1$, we define $\bar{x}_i = [(\bar{x}_i)_1, (\bar{x}_i)_2, \ldots, (\bar{x}_i)_n]$ as follows:

$$(\bar{x}_i)_j = \begin{cases} 1 & a_{ij} \leq b_i^1 \\ \dfrac{(1+\lambda)b_i^1 + 1 - a_{ij}}{1 + \lambda a_{ij}} & a_{ij} > b_i^1 \end{cases}$$

**Proposition 1 [5].** $S_{T_{SW}^\lambda}(A, b^1) = [0, \bar{X}]$ where $0$ is an $n$-dimensioal zero vector and $\bar{X} = \min_{i \in I_1}\{\bar{x}_i\}$.

According to the definition of $S_{T_{SW}^\lambda}(D, b^2)$, the following result is trivially attained.

**Lemma 1.** If $x \in S_{T_{SW}^\lambda}(D, b^2)$, then $x' \in S_{T_{SW}^\lambda}(D, b^2)$ holds for any $x' \in [0,1]^n$ satisfying $x \leq x'$.

The consistency of problem (1) can be checked by the following Proposition.

**Proposition 2.** The following statements are equivalent:

(a) $S_{T_{SW}^\lambda}(D, b^2) \neq \emptyset$. (b) $1 \in S_{T_{SW}^\lambda}(D, b^2)$ where 1 is an $n$-dimensioal vector with each component equal to one. (c) $J_i^2 \neq \emptyset$, $\forall i \in I_2$.

**Proof.** (a)$\Rightarrow$(b) Let $x \in S_{T_{SW}^\lambda}(D, b^2)$. Since $x \leq 1$, then by Lemma 1, $1 \in S_{T_{SW}^\lambda}(D, b^2)$. (b)$\Rightarrow$(c) Suppose that $1 \in S_{T_{SW}^\lambda}(D, b^2)$. So, $max_{j \in J}\{T_{SW}^\lambda(d_{ij}, 1)\} = max_{j \in J}\{d_{ij}\} \geq b_i^2$, $\forall i \in I_2$. Hence, there exists at least one $j_i \in J$ such that $d_{ij_i} \geq b_i^2$. (c)$\Rightarrow$(a) Let $J_i^2 \neq \emptyset$, $\forall i \in I_2$. Thus, for each $i \in I_2$, there exists at least one $j_i \in J$ such that $d_{ij_i} \geq b_i^2$. So, similar to part (a), we have $max_{j \in J}\{T_{SW}^\lambda(d_{ij}, 1)\} = max_{j \in J}\{d_{ij}\} \geq d_{ij_i} \geq b_i^2$, $\forall i \in I_2$. Thus, $1 \in S_{T_{SW}^\lambda}(D, b^2)$. □

**Corollary 1.** Suppose that $S_{T_{SW}^\lambda}(D, b^2) \neq \emptyset$. Then, vector 1 is the maximum solution of $S_{T_{SW}^\lambda}(D, b^2)$.

**3. Characterization of the feasible solution set**

Let $S_{T_{SW}^\lambda}(d_i, b_i^2) = \{x \in [0,1]^n : max_{j \in J}\{T_{SW}^\lambda(d_{ij}, x_j)\} \geq b_i^2\}$, $\forall i \in I_2$, where $d_i$ denotes the $i$th row of matrix $D$. So, it is clear that $S_{T_{SW}^\lambda}(D, b^2) = \bigcap_{i \in I_2} S_{T_{SW}^\lambda}(d_i, b_i^2)$. Based on Corollary 1, vector 1 is the maximum solution of $S_{T_{SW}^\lambda}(d_i, b_i^2)$, $\forall i \in I_2$. The following Proposition shows that the set $S_{T_{SW}^\lambda}(d_i, b_i^2)$ has exactly $|J_i^2|$ minimal solutions, where $|J_i^2|$ denotes the cardinality of the set $J_i^2$.

**Proposition 3.** Suppose that $i \in I_2$, $j_0 \in J_i^2$ and $S_{T_{SW}^\lambda}(d_i, b_i^2) \neq \emptyset$. Also, define $\underline{x}(i, j_0) \in [0,1]^n$ where

$$\underline{x}(i,j_0)_j = \begin{cases} \dfrac{(1+\lambda)b_i^2 + 1 - d_{ij}}{1+\lambda d_{ij}} & j = j_0, b_i^2 > 0 \\ 0 & j = j_0, b_i^2 = 0 \\ 0 & j \neq j_0 \end{cases}$$

Then, $\underline{x}(i,j_0)$ is a minimal solution for $S_{T_{SW}^\lambda}(d_i, b_i^2)$.

**Proof.** By contradiction, suppose that $x' \in S_{T_{SW}^\lambda}(d_i, b_i^2)$, $x' \leq \underline{x}(i,j_0)$ and $x' \neq \underline{x}(i,j_0)$. So, $x'_j \leq \underline{x}(i,j_0)_j$, $\forall j \in J$ and $x' \neq \underline{x}(i,j_0)$. Therefore, $x'_j = 0$, $\forall j \in J - \{j_0\}$. If $x'_j < 0$, then $x' \notin S_{T_{SW}^\lambda}(d_i, b_i^2)$ that is a contradiction. Otherwise, $x'_{j_0} < \frac{(1+\lambda)b_i^2 + 1 - d_{ij_0}}{1+\lambda d_{ij_0}}$. Hence, we have $T_{SW}^\lambda(d_{ij}, 0) = 0$, $\forall j \in J - \{j_0\}$, and $T_{SW}^\lambda(d_{ij_0}, x'_{j_0}) < b_i^2$. Therefore, $max_{j \in J}\{T_{SW}^\lambda(d_{ij}, x'_j)\} < b_i^2$ which contradicts $x' \in S_{T_{SW}^\lambda}(d_i, b_i^2)$. □

**Corollary 2.** $S_{T_{SW}^\lambda}(d_i, b_i^2) = \bigcup_{j \in J_i^2}[\underline{x}(i,j), 1]$, $\forall i \in I_2$.

**Definition 3.** Let $e: I \to \bigcup_{j \in I_2} J_i^2$ so that $e(i) \in J_i^2$, $\forall i \in I_2$, and let $E$ be the set of all vectors $e$. Also, for each $e \in E$, we define $\underline{x}(e) \in [0,1]^n$ such that $\underline{x}(e) = max_{i \in I_2}\{\underline{x}(i, e(i))\}$.

Proposition 4 below shows that the feasible solutions set of $D\varphi x \geq b^2$ can be written in terms of a finite closed convex cells.

**Proposition 4.** $S_{T_{SW}^\lambda}(D, b^2) = \bigcup_{e \in E}[\underline{x}(e), 1]$.

**Proof.** From Corollary 2 and the equality $S_{T_{SW}^\lambda}(D, b^2) = \bigcap_{i \in I_2} S_{T_{SW}^\lambda}(d_i, b_i^2)$, we have $S_{T_{SW}^\lambda}(D, b^2) = \bigcap_{i \in I_2} \bigcup_{j \in J_i^2}[\underline{x}(i,j), 1]$, or equivalently $S_{T_{SW}^\lambda}(D, b^2) = \bigcup_{e \in E} \bigcap_{i \in I_2}[\underline{x}(i, e(i)), 1]$. Therefore, $S_{T_{SW}^\lambda}(D, b^2) = \bigcup_{e \in E}[max_{i \in I_2}\{\underline{x}(i, e(i))\}, 1]$. Now, the result follows from the definition of $\underline{x}(e)$. □

The feasible solutions set of Problem (1) is resulted by the following Proposition.

**Proposition 5.** $S_{T_{SW}^\lambda}(A, D, b^1, b^2) = \bigcup_{e \in E}[\underline{x}(e), \bar{X}]$.

**Proof.** The result follows from Proposition 1, Proposition 4 and the equality $S_{T_{SW}^\lambda}(A, D, b^1, b^2) = S_{T_{SW}^\lambda}(A, b^1) \cap S_{T_{SW}^\lambda}(D, b^2)$. □

**Corollary 3.** Assume that $S_{T_{SW}^\lambda}(A, D, b^1, b^2) \neq \emptyset$. Then, $\bar{X} \in S_{T_{SW}^\lambda}(D, b^2)$.

**4. Optimal solutions of the problem**

Based on Proposition 5, for each $x \in S_{T_{SW}^\lambda}(A, D, b^1, b^2)$ there exist $e' \in E$ such that $x \in [\underline{x}(e'), \bar{X}]$. In other words, if $x \leq \underline{x}(e)$ and $x \neq \underline{x}(e)$, $\forall e \in E$, then $x \notin S_{T_{SW}^\lambda}(A, D, b^1, b^2)$. So, we have the following corollary.

**Corollary 4.** $\underline{S}_{T_{SW}^\lambda}(A, D, b^1, b^2) \subseteq \{\underline{x}(e) : e \in E\}$ where $\underline{S}_{T_{SW}^\lambda}(A, D, b^1, b^2)$ denotes the minimal solutions set of $S_{T_{SW}^\lambda}(A, D, b^1, b^2)$.

The following Proposition characterizes the optimal solution of Problem (1).

**Proposition 6.** If $S_{T_{SW}^\lambda}(A, D, b^1, b^2) \neq \emptyset$, there exists a minimal solution $\underline{x}^* \in \underline{S}_{T_{SW}^\lambda}(A, D, b^1, b^2)$, such that $\underline{x}^*$ is an optimal solution of problem (1).

**Proof.** Let $z(x) = max\{x_1, x_2, \ldots, x_n\}$. Furthermore, suppose that $z(\underline{x}^*) = min\{z(\underline{x}): \underline{x} \in \underline{S}_{T_{SW}^\lambda}(A, D, b^1, b^2)\}$. Based on Proposition 5, for each $x' \in S_{T_{SW}^\lambda}(A, D, b^1, b^2)$ there exist some $\underline{x} \in \underline{S}_{T_{SW}^\lambda}(A, D, b^1, b^2)$ such that $\underline{x} \leq x'$, i.e., $\underline{x}_j \leq x'_j$, $\forall j \in J$. So, we have $\max\{\underline{x}_1, \underline{x}_2, ..., \underline{x}_n\} \leq \max\{x'_1, x'_2, ..., x'_n\}$ that implies $z(\underline{x}) \leq z(x')$. But, $z(\underline{x}^*) \leq z(\underline{x})$ which implies $z(\underline{x}^*) \leq z(x')$, $\forall x' \in S_{T_{SW}^\lambda}(A, D, b^1, b^2)$. □

By combination of Proposition 6 and Corollary 4, it turns out that the optimal solution of problem (1) must be a vector $\underline{x}(e^*)$ for some $e^* \in E$. Based on this fact, we can find the optimal solution of problem (1) by pairwise comparison between the elements of set $\{\underline{x}(e): e \in E\}$. We now summarize the preceding discussion as an algorithm.

**Algorithm 1.** Given Problem (1):

**1.** If $1 \notin S_{T_{SW}^\lambda}(D, b^2)$, then stop; $S_{T_{SW}^\lambda}(A, D, b^1, b^2)$ is empty (Proposition 2).
**2.** If $\bar{X} \notin S_{T_{SW}^\lambda}(D, b^2)$, then stop; $S_{T_{SW}^\lambda}(A, D, b^1, b^2)$ is empty (Corollary 3).
**3.** Find solutions $\underline{x}(e), \forall e \in E$ (Definition 3).
**4.** Find the minimal solutions, $\underline{S}_{T_{SW}^\lambda}(A, D, b^1, b^2)$ by the pairwise comparison between the solutions $\underline{x}(e)$ (Corollary 4).
**5.** Find the optimal solution $\underline{x}(e^*)$ for problem (1) by the pairwise comparison between the objective values of the elements of $\underline{S}_{T_{SW}^\lambda}(A, D, b^1, b^2)$ (Proposition 6).

## 5. Numerical example

Consider the following latticized optimization problem (1):

$$\min \ Z(x) = \max\{x_1, ..., x_{10}\}$$

$$\begin{bmatrix}
0.1576 & 0.6557 & 0.7060 & 0.4387 & 0.2760 & 0.7513 & 0.8407 & 0.3517 & 0.0759 & 0.1622 \\
0.9706 & 0.0357 & 0.0318 & 0.3816 & 0.6797 & 0.2551 & 0.2543 & 0.8308 & 0.0540 & 0.7943 \\
0.9572 & 0.8491 & 0.2769 & 0.7655 & 0.6551 & 0.5060 & 0.8143 & 0.5853 & 0.5308 & 0.3112 \\
0.4854 & 0.9340 & 0.0462 & 0.7952 & 0.1626 & 0.6991 & 0.2435 & 0.5497 & 0.7792 & 0.5285 \\
0.8003 & 0.6787 & 0.0971 & 0.1869 & 0.1190 & 0.8909 & 0.9293 & 0.9172 & 0.9340 & 0.1656 \\
0.1419 & 0.7577 & 0.8235 & 0.4898 & 0.4984 & 0.9593 & 0.3500 & 0.2858 & 0.1299 & 0.6020 \\
0.4218 & 0.7431 & 0.6948 & 0.4456 & 0.9597 & 0.5472 & 0.1966 & 0.7572 & 0.5688 & 0.2630 \\
0.9157 & 0.3922 & 0.3171 & 0.6463 & 0.3404 & 0.1386 & 0.2511 & 0.7537 & 0.4694 & 0.6541 \\
0.7922 & 0.6555 & 0.9502 & 0.7094 & 0.5853 & 0.1493 & 0.6160 & 0.3804 & 0.0119 & 0.6892 \\
0.9595 & 0.1712 & 0.0344 & 0.7547 & 0.2238 & 0.2575 & 0.4733 & 0.5678 & 0.3371 & 0.7482
\end{bmatrix} \varphi x \leq \begin{bmatrix} 0.4505 \\ 0.0838 \\ 0.2290 \\ 0.9133 \\ 0.1524 \\ 0.8258 \\ 0.5383 \\ 0.9961 \\ 0.0782 \\ 0.4427 \end{bmatrix}$$

$$\begin{bmatrix}
0.3674 & 0.4899 & 0.9730 & 0.6099 & 0.3909 & 0.4588 & 0.6797 & 0.9810 & 0.8909 & 0.0714 \\
0.9880 & 0.1679 & 0.6490 & 0.6177 & 0.8314 & 0.7984 & 0.1366 & 0.0855 & 0.3342 & 0.5216 \\
0.9995 & 0.9787 & 0.8003 & 0.8594 & 0.8034 & 0.9631 & 0.7212 & 0.2625 & 0.6987 & 0.0967 \\
0.0377 & 0.7127 & 0.4538 & 0.8055 & 0.0605 & 0.5468 & 0.7982 & 0.8010 & 0.1978 & 0.8181 \\
0.8852 & 0.5005 & 0.4324 & 0.5767 & 0.3993 & 0.5211 & 0.1068 & 0.0292 & 0.8241 & 0.8175 \\
0.9133 & 0.4711 & 0.8253 & 0.8643 & 0.5269 & 0.2316 & 0.6538 & 0.9289 & 0.0305 & 0.7224 \\
0.7962 & 0.9426 & 0.0835 & 0.1829 & 0.4168 & 0.4889 & 0.4942 & 0.7303 & 0.7441 & 0.1499 \\
0.0987 & 0.0596 & 0.1332 & 0.2399 & 0.6569 & 0.6241 & 0.7791 & 0.4886 & 0.5000 & 0.7311 \\
0.2619 & 0.6820 & 0.8822 & 0.8865 & 0.6280 & 0.6791 & 0.7150 & 0.5785 & 0.4799 & 0.6596 \\
0.3354 & 0.0424 & 0.1734 & 0.0287 & 0.8612 & 0.3955 & 0.9037 & 0.2373 & 0.9047 & 0.5186
\end{bmatrix} \varphi x \geq \begin{bmatrix} 0.0492 \\ 0.0876 \\ 0.0864 \\ 0.0810 \\ 0.0471 \\ 0.1292 \\ 0.1511 \\ 0.0392 \\ 0.0291 \\ 0.1028 \end{bmatrix}$$

$x \in [0,1]^{10}$

where $|I_1| = |I_2| = |J| = 10$ and $\varphi(x, y) = T_{SW}^\lambda(x, y)$ in which $\lambda = 2$.

In this example, Also, from Definition 1, $J_1^2 = J_3^2 = J_8^2 = J_9^2 = J = \{1, ..., 10\}$, $J_2^2 = J_5^2 = J - \{8\}$, $J_6^2 = J - \{9\}$, $J_7^2 = J - \{3,10\}$ and $J_{10}^2 = J - \{2,4\}$. By Definition 2, we have

$\bar{x}_1 = [1, 0.7337, 0.6822, 1, 1, 0.6395, 0.5635, 1, 1, 1]$,
$\bar{x}_2 = [0..955, 1, 1, 0.4934, 0.2423, 0.6598, 0.661, 0.158, 1, 0.1766]$,
$\bar{x}_3 = [0.2504, 0.3105, 0.9074, 0.364, 0.4466, 0.587, 0.332, 0.5076, 0.5608, 0.8479]$,
$\bar{x}_4 = [1, 0.9784, 1, 1, 1, 1, 1, 1, 1, 1]$,
$\bar{x}_5 = [0.2526, 0.3302, 1, 0.9247, 1, 0.2035, 0.1847, 0.1905, 0.1824, 0.9701]$

$\bar{x}_6 = [1,1,1,1,1,0.8628,1,1,1,1]$
$\bar{x}_7 = [1,0.7529,0.8035,1,0.567,0.9873,1,0.7389,0.9572,1]$
$\bar{x}_8 = [1,1,1,1,1,1,1,1,1,1]$
$\bar{x}_9 = [0.1711,0.2506,0.098,0.2171,0.2991,0.8357,0.2771,0.485,1,0.2293]$
$\bar{x}_{10} = [0.4688,1,1,0.627,1,1,0.9528,0.8242,1,0.6329]$

From Proposition 2, the necessary condition holds for the feasibility of the problem. More precisely, we have $D\varphi 1 \geq b^2$ that means $1 \in S_{T_{SW}^\lambda}(D, b^2)$.

By Proposition 1, $\bar{X} = [0.4496,0.0949,0.4138,0.1918,0.1991,0.1731,0.8846,0,0.1437,0]$ which determines the feasible region of the first inequalities, i.e., $S_{T_{SW}^\lambda}(A, b^1) = [0, \bar{X}]$. Also, $D\varphi\bar{X} \geq b^2$.

Therefore, we have $\bar{X} \in S_{T_{SW}^\lambda}(D, b^2)$, which satisfies the necessary feasibility condition stated in Corollary 3. On the other hand, from Definition 3, we have $|E| = 3{,}732{,}480{,}000$. Therefore, the number of all vectors $e \in E$ is equal to 3,732,480,000. In order to find the optimal solution $\underline{x}(e^*)$ of sub-problems (4), we firstly compute all minimal solutions by making pairwise comparisons between all solutions $\underline{x}(e)$ ($\forall e \in E$), and then we find $\underline{x}(e^*)$ among the resulted minimal solutions. Actually, the feasible region has 38 minimal solutions reported in Appendix section. By comparison of the values of the objective function for the minimal solutions, we have $z(\underline{x}(e_i)) = 0.1991$ for $i \in \{6,8,16,19,21,29,33\}$. Hence, this problem has 7 optimal solutions $\underline{x}(e_i)$ for $i \in \{6,8,16,19,21,29,33\}$.

**Conclusion**

Considering the practical applications of the max-Sugeno-Weber fuzzy relational inequalities in FRI theory and that of the latticized problems in the convex optimization, a nonlinear optimization problem was studied with the maximum components function as the objective function subjected to the Sugeno-Weber-FRI. Since a system of the Sugeno-Weber-FRI is a non-convex set, an algorithm was presented to find an optimal solution by using the structural properties of the problem. For this purpose, some feasibility conditions were firstly derived and then, the feasible region was completely determined in terms of one maximum and a finite number of minimal solutions. It is proved that we can find the exact optimal solution of the proposed problem from the minimal solutions of the constraints, i.e., a system of max-Sugeno-Weber FRI. Additionally, a numerical example was given to illustrate the presented algorithm.

**Appendix**

The minimal solutions of the feasible solutions set for the problem stated in section 5.

$e_1 = [3, 1, 1, 4, 1, 4, 2, 7, 2, 7]$
$\underline{x}(e_1) = [0.0924 \quad 0.1770 \quad 0.0592 \quad 0.1918 \quad 0 \quad 0 \quad 0.1441 \quad 0 \quad 0 \quad 0]$

$e_2 = [8, 1, 1, 4, 1, 4, 2, 7, 2, 7]$
$\underline{x}(e_2) = [0.0924 \quad 0.1770 \quad 0 \quad 0.1918 \quad 0 \quad 0 \quad 0.1441 \quad 0.0562 \quad 0 \quad 0]$

$e_3 = [3, 1, 1, 4, 1, 4, 2, 7, 2, 5]$
$\underline{x}(e_3) = [0.0924 \quad 0.1770 \quad 0.0592 \quad 0.1918 \quad 0.1642 \quad 0 \quad 0.1323 \quad 0 \quad 0 \quad 0]$

$e_4 = [8, 1, 1, 4, 1, 4, 2, 7, 2, 5]$
$\underline{x}(e_4) = [0.0924 \quad 0.1770 \quad 0 \quad 0.1918 \quad 0.1642 \quad 0 \quad 0.1323 \quad 0.0562 \quad 0 \quad 0]$

$e_5 = [9, 6, 2, 4, 9, 4, 2, 7, 2, 7]$
$\underline{x}(e_5) = [0 \quad 0.1770 \quad 0 \quad 0.1918 \quad 0 \quad 0.1788 \quad 0.1441 \quad 0 \quad 0.1197 \quad 0]$

$e_6 = [3, 1, 1, 4, 1, 4, 2, 5, 2, 5]$
$\underline{x}(e_6) = [0.0924 \quad 0.1770 \quad 0.0592 \quad 0.1918 \quad 0.1991 \quad 0 \quad 0 \quad 0 \quad 0 \quad 0]$

$e_7 = [3, 6, 2, 4, 10, 4, 2, 7, 2, 7]$
$\underline{x}(e_7) = [0 \quad 0.1770 \quad 0.0592 \quad 0.1918 \quad 0 \quad 0.1788 \quad 0.1441 \quad 0 \quad 0 \quad 0.1228]$

$e_8 = [8, 1, 1, 4, 1, 4, 2, 5, 2, 5]$
$\underline{x}(e_8) = [0.092 \quad 0.1770 \quad 0 \quad 0.1918 \quad 0.1991 \quad 0 \quad 0 \quad 0.0562 \quad 0 \quad 0]$

$e_9 = [8, 6, 2, 4, 10, 4, 2, 7, 2, 7]$
$\underline{x}(e_9) = [0 \quad 0.1770 \quad 0 \quad 0.1918 \quad 0 \quad 0.1788 \quad 0.1441 \quad 0.0562 \quad 0 \quad 0.1228]$

$e_{10} = [9, 6, 2, 4, 9, 4, 2, 7, 2, 9]$
$\underline{x}(e_{10}) = [0 \quad 0.1770 \quad 0 \quad 0.1918 \quad 0 \quad 0.1788 \quad 0.1323 \quad 0 \quad 0.1437 \quad 0]$

$e_{11} = [9, 1, 1, 4, 1, 4, 2, 7, 2, 7]$
$\underline{x}(e_{11}) = [0.0924 \quad 0.1770 \quad 0 \quad 0.1918 \quad 0 \quad 0 \quad 0.1441 \quad 0 \quad 0.0922 \quad 0]$

$e_{12} = [9, 5, 2, 4, 9, 4, 2, 7, 2, 7]$
$\underline{x}(e_{12}) = [0 \quad 0.1770 \quad 0 \quad 0.1918 \quad 0.1620 \quad 0 \quad 0.1441 \quad 0 \quad 0.1197 \quad 0]$

$e_{13} = [9, 5, 2, 4, 9, 4, 2, 7, 2, 5]$
$\underline{x}(e_{13}) = [0 \quad 0.1770 \quad 0 \quad 0.1918 \quad 0.1642 \quad 0 \quad 0.1323 \quad 0 \quad 0.1197 \quad 0]$

$e_{14} = [3, 5, 2, 4, 10, 4, 2, 7, 2, 7]$
$\underline{x}(e_{14}) = [0 \quad 0.1770 \quad 0.0592 \quad 0.1918 \quad 0.1620 \quad 0 \quad 0.1441 \quad 0 \quad 0 \quad 0.1228]$

$e_{15} = [3, 5, 2, 4, 10, 4, 2, 7, 2, 5]$
$\underline{x}(e_{15}) = [0 \quad 0.1770 \quad 0.0592 \quad 0.1918 \quad 0.1642 \quad 0 \quad 0.1323 \quad 0 \quad 0 \quad 0.1228]$

$e_{16} = [9, 5, 2, 4, 9, 4, 2, 5, 2, 5]$
$\underline{x}(e_{16}) = [0 \quad 0.1770 \quad 0 \quad 0.1918 \quad 0.1991 \quad 0 \quad 0 \quad 0 \quad 0.1197 \quad 0]$

$e_{17} = [8, 5, 2, 4, 10, 4, 2, 7, 2, 7]$
$\underline{x}(e_{17}) = [0 \quad 0.1770 \quad 0 \quad 0.1918 \quad 0.1620 \quad 0 \quad 0.1441 \quad 0.0562 \quad 0 \quad 0.1228]$

$e_{18} = [8, 5, 2, 4, 10, 4, 2, 7, 2, 5]$
$\underline{x}(e_{18}) = [0 \quad 0.1770 \quad 0 \quad 0.1918 \quad 0.1642 \quad 0 \quad 0.1323 \quad 0.0562 \quad 0 \quad 0.1228]$

$e_{19} = [3, 5, 2, 4, 10, 4, 2, 5, 2, 5]$
$\underline{x}(e_{19}) = [0 \quad 0.1770 \quad 0.0592 \quad 0.1918 \quad 0.1991 \quad 0 \quad 0 \quad 0 \quad 0 \quad 0.1228]$

$e_{20} = [9, 1, 1, 4, 9, 4, 2, 7, 2, 7]$
$\underline{x}(e_{20}) = [0.0923 \quad 0.1770 \quad 0 \quad 0.1918 \quad 0 \quad 0 \quad 0.1441 \quad 0 \quad 0.1197 \quad 0]$

$e_{21} = [8, 5, 2, 4, 10, 4, 2, 5, 2, 5]$
$\underline{x}(e_{21}) = [0 \quad 0.1770 \quad 0 \quad 0.1918 \quad 0.1991 \quad 0 \quad 0 \quad 0.0562 \quad 0 \quad 0.1228]$

$e_{22} = [3, 1, 1, 4, 10, 4, 2, 7, 2, 7]$
$\underline{x}(e_{22}) = [0.0923 \quad 0.1770 \quad 0.0592 \quad 0.1918 \quad 0 \quad 0 \quad 0.1441 \quad 0 \quad 0 \quad 0.1228]$

$e_{23} = [9, 5, 2, 4, 9, 4, 2, 7, 2, 9]$
$\underline{x}(e_{23}) = [0 \quad 0.1770 \quad 0 \quad 0.1918 \quad 0.1620 \quad 0 \quad 0.1323 \quad 0 \quad 0.1437 \quad 0]$

$e_{24} = [8, 1, 1, 4, 10, 4, 2, 7, 2, 7]$
$\underline{x}(e_{24}) = [0.0923 \quad 0.1770 \quad 0 \quad 0.1918 \quad 0 \quad 0 \quad 0.1441 \quad 0.0562 \quad 0 \quad 0.1228]$

$e_{25} = [9, 1, 1, 4, 1, 4, 2, 7, 2, 5]$
$\underline{x}(e_{25}) = [0.0924 \quad 0.1770 \quad 0 \quad 0.1918 \quad 0.1642 \quad 0 \quad 0.1323 \quad 0 \quad 0.0922 \quad 0]$

$e_{26} = [3, 5, 2, 4, 10, 4, 2, 10, 2, 5]$
$\underline{x}(e_{26}) = [0 \quad 0.1770 \quad 0.0592 \quad 0.1918 \quad 0.1642 \quad 0 \quad 0 \quad 0 \quad 0 \quad 0.1570]$

$e_{27} = [9, 1, 1, 4, 9, 4, 2, 7, 2, 9]$
$\underline{x}(e_{27}) = [0.0923 \quad 0.1770 \quad 0 \quad 0.1918 \quad 0 \quad 0 \quad 0.1323 \quad 0 \quad 0.1437 \quad 0]$

$e_{28} = [8, 5, 2, 4, 10, 4, 2, 10, 2, 5]$
$\underline{x}(e_{28}) = [0 \quad 0.1770 \quad 0 \quad 0.1918 \quad 0.1642 \quad 0 \quad 0 \quad 0.0562 \quad 0 \quad 0.1570]$

$e_{29} = [9, 1, 1, 4, 1, 4, 2, 5, 2, 5]$
$\underline{x}(e_{29}) = [0.0924 \quad 0.1770 \quad 0 \quad 0.1918 \quad 0.1991 \quad 0 \quad 0 \quad 0 \quad 0.0922 \quad 0]$

$e_{30} = [9, 6, 2, 4, 10, 4, 2, 7, 2, 7]$
$\underline{x}(e_{30}) = [0 \quad 0.1770 \quad 0 \quad 0.1918 \quad 0 \quad 0.1788 \quad 0.1441 \quad 0 \quad 0.0922 \quad 0.1228]$

$e_{31} = [9, 5, 2, 4, 10, 4, 2, 7, 2, 7]$
$\underline{x}(e_{31}) = [0 \quad 0.1770 \quad 0 \quad 0.1918 \quad 0.1620 \quad 0 \quad 0.1441 \quad 0 \quad 0.0922 \quad 0.1228]$

$e_{32} = [9, 5, 2, 4, 10, 4, 2, 7, 2, 5]$
$\underline{x}(e_{32}) = [0 \quad 0.1770 \quad 0 \quad 0.1918 \quad 0.1642 \quad 0 \quad 0.1323 \quad 0 \quad 0.0922 \quad 0.1228]$

$e_{33} = [9, 5, 2, 4, 10, 4, 2, 5, 2, 5]$
$\underline{x}(e_{33}) = [0 \quad 0.1770 \quad 0 \quad 0.1918 \quad 0.1991 \quad 0 \quad 0 \quad 0 \quad 0.0922 \quad 0.1228]$

$e_{34} = [9, 1, 1, 4, 10, 4, 2, 7, 2, 7]$
$\underline{x}(e_{34}) = [0.0923 \quad 0.1770 \quad 0 \quad 0.1918 \quad 0 \quad 0 \quad 0.1441 \quad 0 \quad 0.0922 \quad 0.1228]$

$e_{35} = [9, 5, 2, 4, 10, 4, 2, 10, 2, 5]$
$\underline{x}(e_{35}) = [0 \quad 0.1770 \quad 0 \quad 0.1918 \quad 0.1642 \quad 0 \quad 0 \quad 0 \quad 0.0922 \quad 0.1570]$

$e_{36} = [9, 6, 2, 4, 9, 4, 2, 10, 2, 9]$
$\underline{x}(e_{36}) = [0 \quad 0.1770 \quad 0 \quad 0.1918 \quad 0 \quad 0.1788 \quad 0 \quad 0 \quad 0.1437 \quad 0.1570]$

$e_{37} = [9, 5, 2, 4, 9, 4, 2, 10, 2, 9]$
$\underline{x}(e_{37}) = [0 \quad 0.1770 \quad 0 \quad 0.1918 \quad 0.1620 \quad 0 \quad 0 \quad 0 \quad 0.1437 \quad 0.1570]$

$e_{38} = [9, 1, 1, 4, 9, 4, 2, 10, 2, 9]$
$\underline{x}(e_{38}) = [0.0923 \quad 0.1770 \quad 0 \quad 0.1918 \quad 0 \quad 0 \quad 0 \quad 0 \quad 0.1437 \quad 0.1570]$

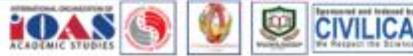